\documentclass[a4paper,12pt]{article}
\usepackage{latexsym,amsmath,amssymb,amscd}
\begin{document}
\title{Fourier Algebra of a Compact Lie Group}
\author{R.J. Plymen}
\date{}
\maketitle

Let $G$ be a compact group and let $A(G)$ be the Fourier algebra
of $G$.  This is a commutative Banach algebra.   As a Banach space,
$A(G)$ is defined as the predual of $vN(G)$, the von Neumann algebra
of $G$.  Now weak amenability of $A(G)$ in the sense of \cite{BCD} has
been studied by Johnson \cite{J}, and we strengthen one of his results
with the following result.

\textsc{Theorem}.  Let $G$ be a compact connected
Lie group.  Then $A(G)$ is weakly amenable if and only if $G$ is
abelian.
                 
\textsc{Proof} (i)  Let $G$ be abelian.  Then the Pontryagin dual
$\hat G$ is a discrete abelian group.   Then we have
$vN(G) \cong \ell^{\infty}(\hat G)$ and so we have

\[A(G) \cong \ell^{1} (\hat G)\]

Now $\ell^{1} (\hat G)$ is amenable.   Therefore
$A(G)$ is weakly amenable.

(ii)  Let $G$ be a non-abelian compact connected Lie group.
For Lie groups our reference is Bourbaki \cite{B} especially
Lie IX.5 and Lie IX.31.  The Lie algebra $L(G)$ decomposes as

\[L(G) = c \oplus D(L(G))\]

where $c$ is the centre of $L(G)$ and $D(L(G))$ is the
derived algebra of $L(G)$.  The derived algebra $D(L(G))$ is a
semisimple Lie algebra.  If $G$ is non-abelian then $L(G)$ is
non-abelian and $D(L(G)) \neq 0$.

Hence the root system

\[R(g_\mathbb{C}, t_\mathbb{C}) \neq 0\]
where $g = L(G), t = L(T)$ and $T$ is a maximal torus of $G$.
But $R(G,T) \cong R(g_\mathbb{C}, t_\mathbb{C})$.  Hence there exists a
root $\alpha \in R(G,T)$.

Let $Z_{\alpha}$ be the centralizer of the kernel of $\alpha$.
Let $S_{\alpha} = D(Z_{\alpha})$ be the derived group (commutator
subgroup) of $Z_{\alpha}$.  Then $S_{\alpha}$ is a connected compact
semisimple Lie group of rank 1.

Let $\pi$ be an irreducible representation of $G$ of degree $d$:
           \[\pi: G \longrightarrow U(d)\]

Consider $T \cap S_{\alpha}$.  This is a maximal torus of
$S_{\alpha}$.   There are only 2 possibilities for $S_{\alpha}$.

(i) $S_{\alpha} = SU(2,\mathbb{C})$.  Consider the restriction
$\pi |SU(2,\mathbb{C})$ and the fact that $T \cap S_{\alpha}$ is a
circle $U(1)$.  Now

\[\pi |SU(2) = \lambda_{1} \oplus \ldots \oplus \lambda_{r}\]

with each $\lambda_{1}, \ldots, \lambda_{r}$ irreducible.
But $\lambda_{j}|U(1) = \chi_{-m} + \ldots + \chi_{m}$ 
where $\chi_{m} (z) = z^{m}$.  Consequently we have $|m| < d$.

(ii) $S_{\alpha} = SO(3,\mathbb{R})$.  A similar argument gives
$|m| < d$ for the characters occurring in $\pi |U(1)$.

We now apply the remark in Johnson \cite{J}(p.373) to infer
that $A(G)$ is not weakly amenable.

\textsc{Corollary}.  If $G$ is any compact
connected Lie group except $SO(2)$ then $A(G)$ is not
weakly amenable.

Mathematics Department, University of Manchester, Manchester M13 9PL, UK.
                  

\begin{thebibliography}{99}
\bibitem{BCD}  W.G. Bade, P.C. Curtis, Jr and H.G. Dales.
Amenability and weak amenability for Buerling and Lipschitz
algebras, Proc. London Math.Soc. 55 (1987) 359 - 377.
\bibitem{B}  N. Bourbaki.  Groupes et alg\`{e}bres de Lie.
Chapitre 9.  Masson, Paris 1982.
\bibitem{J}   B.E. Johnson.  Non-amenability of the Fourier
algebra of a compact group.   J. London Math. Soc. 50 (1994) 361 --374.
\end{thebibliography}
\end{document}